\documentclass[12pt,letterpaper]{amsart}

\usepackage{hyperref}
\usepackage{palatino, euler, epic,eepic, amssymb, xypic, floatflt, microtype}
\usepackage{verbatim,color}
\usepackage{mathtools}
\usepackage{amsmath}
\usepackage{tikz-cd}
\usepackage{listings}
\usepackage{url}

\input xy
\xyoption{all}

\newcommand{\bpf}{\noindent {\em Proof.  }}
\newcommand{\epf}{\qed \vspace{+10pt}}

\newcommand{\tpoint}[1]{\vspace{3mm}\par \noindent \refstepcounter{subsection}{\thesubsection.} 
  {\bf #1. ---} }
\newcommand{\epoint}[1]{\vspace{3mm}\par \noindent \refstepcounter{subsection}{\thesubsection.} 
  {\em #1.} }
\newcommand{\bpoint}[1]{\vspace{3mm}\par \noindent \refstepcounter{subsection}{\thesubsection.} 
  {\bf #1.} }

\def\CC{\mathbb{C}}

\def\QQ{\mathbb{Q}}
\def\RR{\mathbb{R}}

\def\ZZ{\mathbb{Z}}

\newcommand{\cF}{{\mathscr{F}}}

\renewcommand{\cL}{{\mathscr{L}}}

\newcommand{\cS}{{\mathcal{S}}}

\newcommand{\propernormal}{%
  \mathrel{\ooalign{$\lneq$\cr\raise.22ex\hbox{$\lhd$}\cr}}}

\newcommand{\Gal}{\operatorname{Gal}}

\newcommand{\lremind}[1]{{\bf[label:  #1]}}
\newcommand{\notation}[1]{}
\renewcommand{\lremind}[1]{{}}

\newcommand{\cut}[1]{}

\DeclareMathOperator{\GL}{GL}
\DeclareMathOperator{\GO}{GO}

\DeclareMathOperator{\Reg}{Reg}

\DeclareMathOperator{\Log}{Log}
\DeclareMathOperator{\gl}{GL}
\DeclareMathOperator{\mfh}{\mathfrak{h}}

\DeclareMathOperator{\ord}{ord}

\DeclareMathOperator{\Nm}{Nm}

\DeclareMathOperator{\Disc}{Disc}

\DeclareMathOperator{\Covol}{Covol}

\newcommand{\cO}{\mathcal{O}}      

\DeclareSymbolFont{extraup}{U}{zavm}{m}{n}
\DeclareMathSymbol{\varheart}{\mathalpha}{extraup}{86}
\DeclareMathSymbol{\vardiamond}{\mathalpha}{extraup}{87}

\begin{document}
\pagestyle{plain}
\title{\Large{Unit lattices of $D_4$-quartic number fields with signature $(2,1)$}}

\author{Sara Chari}
\address{Dept. of Mathematics, St. Mary’s College of Maryland, Mary’s City, MD~20653}
\email{slchari@smcm.edu}

\author{Sergio Ricardo Zapata Ceballos}
\address{Dept. of Mathematics and Statistics, Youngstown State Uiversity, Youngstown, OH~44555}
\email{srzapataceballos@ysu.edu}

\author{Erik Holmes}
\address{Dept. of Mathematics, University of Toronto, Toronto, ON~M5S 2E4}
\email{eholmes@math.toronto.edu}

\author{Fatemeh Jalalvand}
\address{Dept. of Mathematics, University of Calgary, Calgary, AB~T2N 1N4}
\email{fatemeh.jalalvand@ucalgary.ca}

\author{Rahinatou Yuh Njah Nchiwo}
\address{Dept. of Mathematics and Systems Analysis, Aalto University, Espoo, Finland~02150}
\email{rahinatou.njahepousenchiwo@aalto.fi}

\author{Kelly O'Connor}
\address{Dept. of Mathematics and Statistics, University of Wisconsin-La Crosse, La Crosse WI~54601}
\email{koconnor@uwlax.edu}

\author{Fabian Ramirez}
\address{Dept. of Mathematics, University of California Irvine, Irvine, CA~92697}
\email{fabiar3@uci.edu}

\author{Sameera Vemulapalli}
\address{Dept. of Mathematics, Harvard University, Cambridge MA~02138}
\email{vemulapalli@math.harvard.edu}

\date{\today}

\begin{abstract}
There has been a recent surge of interest on distributions of shapes of unit lattices in number fields, due to both their applications to number theory and the lack of known results. 

In this work we focus on $D_4$-quartic fields with signature $(2,1)$; such fields have a rank $2$ unit group. Viewing the unit lattice as a point of $\gl_2(\ZZ)\backslash \mfh$, we prove that every lattice which arises this way must correspond to a transcendental point on the boundary of a certain fundamental domain of $\gl_2(\ZZ)\backslash \mfh$. Moreover, we produce three explicit (algebraic) points of $\gl_2(\ZZ)\backslash \mfh$ which are limit points of the set of (points associated to) unit lattices of $D_4$-quartic fields with signature $(2,1)$. 
 \end{abstract}
 
\maketitle

\begin{figure}[ht]
\label{fig:1}
\centering
\includegraphics[scale=0.45]{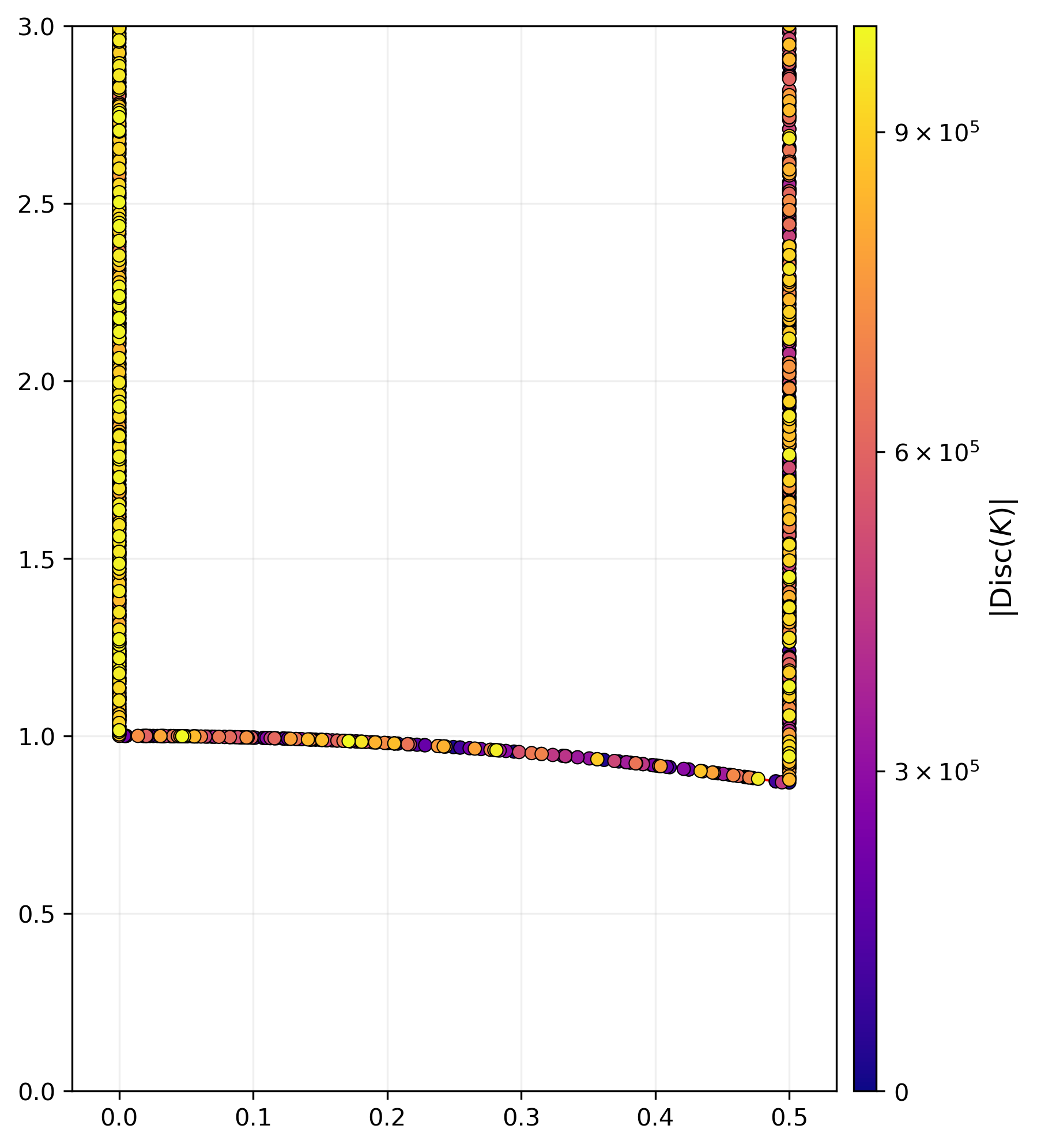}
\caption{Each plotted point corresponds to the shape of the unit lattice of a $D_4$-quartic field with signature $(2,1)$. We list all such fields of absolute discriminant $\leq 10^6$. The points appear to be dense in the boundary of the fundamental domain of $\GL_2(\ZZ) \backslash \mfh$.}\end{figure}

\tableofcontents

{\parskip=12pt 

\section{Introduction}
\label{sec:intro}

Our interest in the study of distributions of shapes of unit lattices in number fields began with the MathOverflow post \cite{mo}, in which A. Reznikov poses a very fundamental set of questions (communicated to him by M. Gromov) on the shape of unit lattices in number fields. 

Namely, we are interested in the following question. For a permutation group $G \subseteq S_n$ and a pair $(r,s)$ of positive integers with $r + 2s = n$, let $\cF(G,r,s)$ be the set of (isomorphism classes) of degree $n$ number fields with Galois group $G$ and signature $(r,s)$ (by standard abuse of notation in arithmetic statistics, we say a number field has Galois group $G$ if its Galois closure has Galois group $G$). For any field $L \in \cF(G,r,s)$, let $\sigma_1,\dots,\sigma_r$ be the real embeddings of $L$ and $\tau_1, \overline{\tau_1},\dots, \tau_s, \overline{\tau_s}$ be the pairs of complex conjugate embeddings of $L$. Let $E_L$ be the group of units of $\cO_K$ modulo torsion. There is a \emph{logarithmic embedding}
\[
    \Log \colon E_L		  	\longrightarrow \RR^{r + s}	
\]
\[
u		\mapsto (\log\lvert \sigma_1(u) \rvert, \hdots, \log\lvert\sigma_r(u)\rvert, 2\log \lvert \tau_1(u) \rvert,\hdots,2\log \lvert \tau_s(u) \rvert).
\]
Restricting the standard quadratic form on $\RR^{r+s}$ to $E_L$ equips $E_L$ with the structure of a lattice of rank $r + s - 1$. We say $E_L$ is the \emph{unit lattice} of $L$. 

We define the \emph{unit shape} of $L$ to be the equivalence class of this lattice under scaling, reflection, and rotation.  Hence for each such number field $L$ we obtain a point, which we call $p_L$, in the moduli space $\mathcal{S}_{r+s-1} = \GL_{r+ s-1}(\ZZ)\backslash \GL_{r+ s-1}(\RR) / \GO_{r+ s-1}(\RR)$ of unimodular rank $r + s - 1$ lattices up to homothety and reflection. We call $p_L$ the \emph{unit shape} of $L$. Define the set:
\[
    \Omega(G,r,s) \coloneqq \{p_L : L \in \cF(G,r,s)\}.
\]

A very natural question is: given a group $G$ and a signature $(r,s)$, where do the unit shapes of $G$-extensions with signature $(r,s)$ lie in $\mathcal{S}_{r+s-1}$? This is, a priori, a very difficult question -- for example, proving that $\Omega(G,r,s) \neq \emptyset$ implies that there exists a $G$-extension of signature $(r,s)$. In particular, understanding when $\Omega(G,r,s) \neq \emptyset$ would provide an answer to the inverse Galois problem. 

\tpoint{Main Question} { \em
\label{question:main-q}
What is $\Omega(G,r,s)$?
}

One could, of course, ask about how the points $p_L$ are distributed as $\cF(G,r,s)$ is ordered in different ways (e.g., by discriminant, regulator, or conductor), but for simplicity we restrict ourselves to the easier question of what the set  $\Omega(G,r,s)$ is. We plot three examples: Figure $1$ is a plot of $\Omega(D_4,2,1)$, Figure $2$ is a plot of $\Omega(S_3,3,0)$, and  Figure $3$ is a plot of the subset of $\Omega(D_4,2,1)$ corresponding to $D_4$-quartic fields containing the fixed real quadratic subfield $\QQ(\sqrt{23})$. We note that, in each of these plots, the color of the point represents the magnitude of the absolute discriminant, $\lvert \Disc(K)\rvert$.

\begin{figure}
\label{fig:combined}
\centering
\begin{minipage}{.45\textwidth}
\label{fig:2}
  \centering
  \includegraphics[width=.9\linewidth]{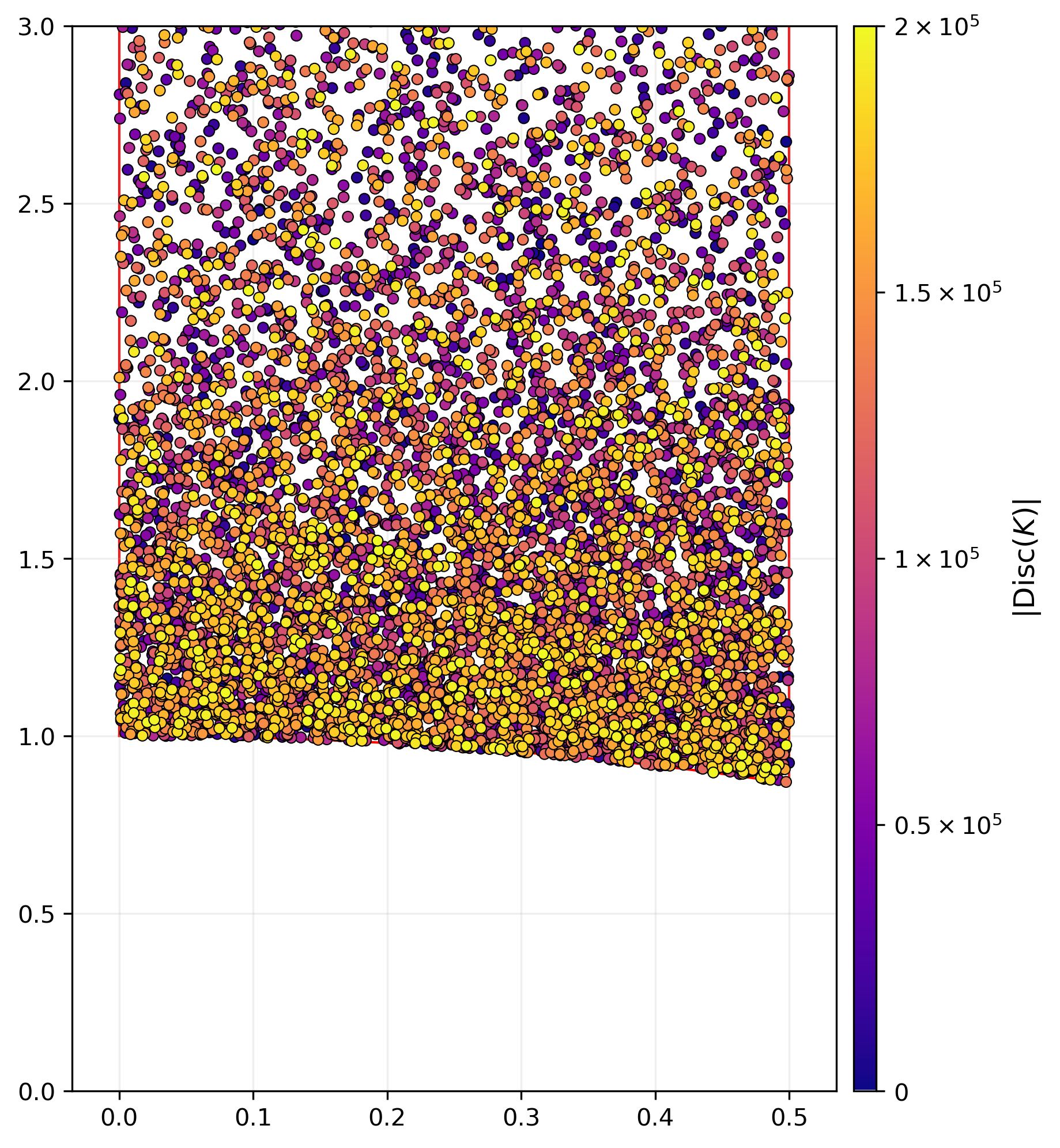}
  \caption{}{The shapes of unit lattices of totally real $S_3$-cubic fields with absolute discriminant $\leq 2\times10^5$. The points appear to be dense in the fundamental domain of $\GL_2(\ZZ) \backslash \mfh$.}
  \label{fig:test1}
\end{minipage}%
\hspace{0.05\textwidth}
\begin{minipage}{.45\textwidth}
\label{fig:3}
  \centering
  \includegraphics[width=.9\linewidth]{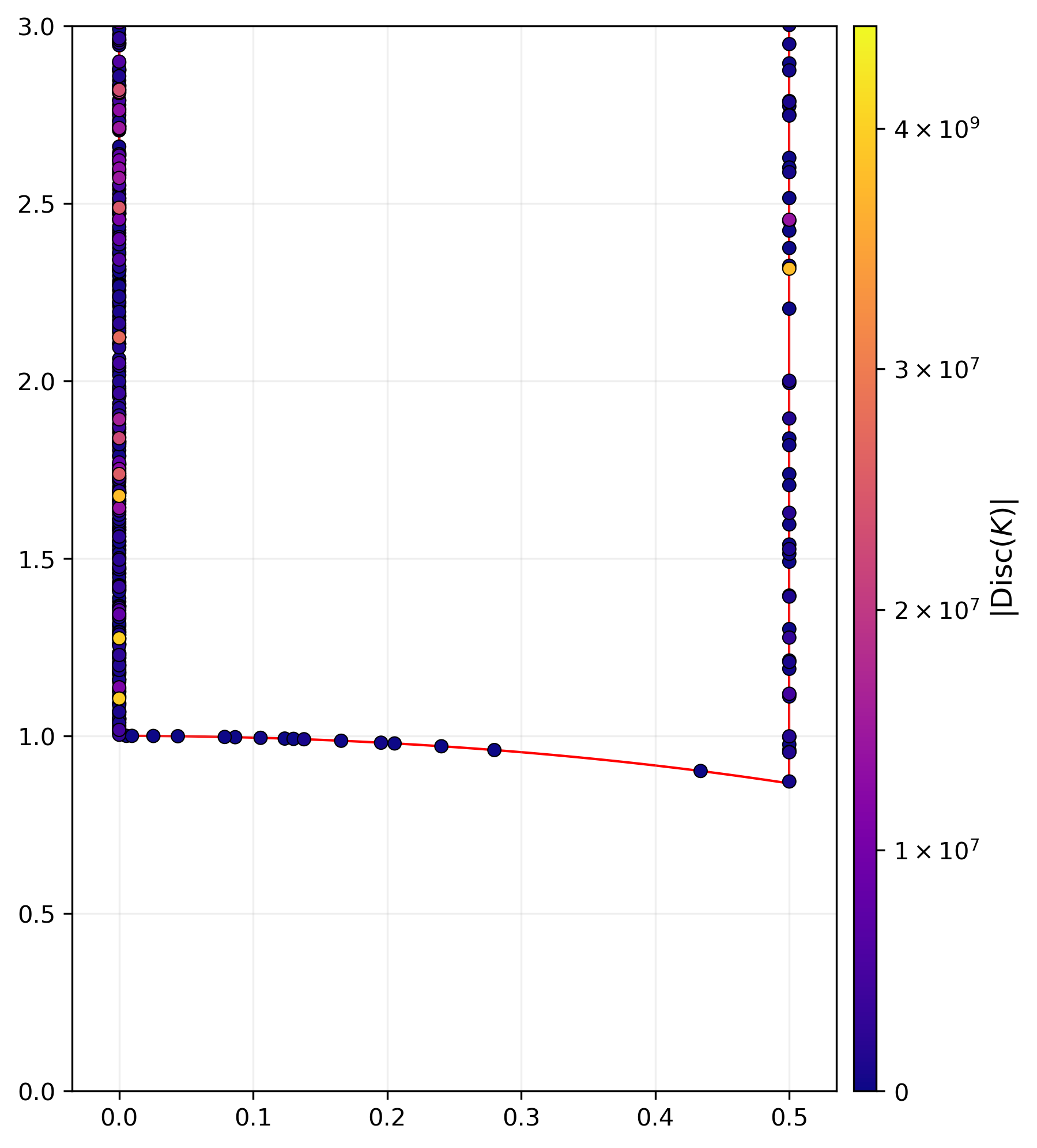}
  \caption{}{The shapes of unit lattices of $D_4$-quartic fields with signature $(2,1)$ and absolute discriminant $\leq 4.5\times10^9$ containing the real quadratic subfield $\mathbb{Q}(\sqrt{23})$. The points appear to be going to infinity.}
  \label{fig:test2}
\end{minipage}
\end{figure}

\bpoint{Previous Work} {
\label{subsec:prev-work}

Very little is known about Question~\ref{question:main-q}. The question is trivial in quadratic fields, as the unit rank is $\leq 1$. So, the first case to consider is the case of unit rank $2$; in rank $2$, the the first signature and degree to consider is that of totally real cubic fields. In this case, we identify $\cS_2$ with the standard fundamental domain, which we call $\cF$, of $\GL_2(\ZZ) \backslash \mfh$ (see Figure $4$). 

\begin{figure}[ht]
\label{fig:4}
\centering
\includegraphics[scale=0.45]{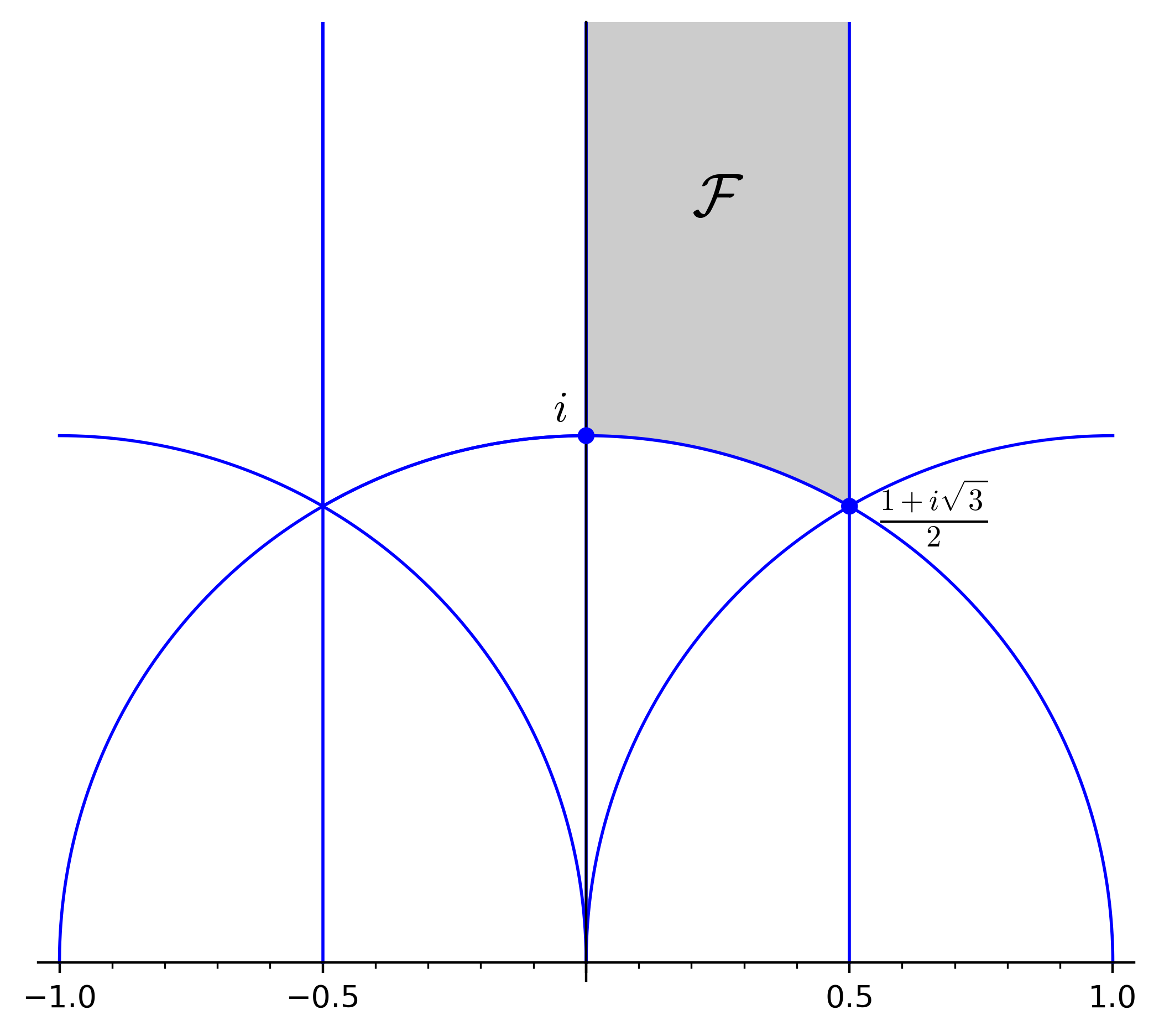}
\caption{The gray region is the standard fundamental domain $\cF$ of $\GL_2(\ZZ) \backslash \mfh$ in the complex plane.}
\end{figure}

When $G = \ZZ/3\ZZ$, the unit lattice has a nontrivial action of $G$, and hence must be the hexagonal lattice, whose corresponding point in $\cF$ is the rightmost cusp $\frac{1}{2} + \frac{i\sqrt{3}}{2}$. In other words, upon identifying $\cS_2$ with the standard fundamental domain of $\GL_2(\ZZ) \backslash \mfh$, we have
\[
    \Omega(\ZZ/3\ZZ,3,0) = \frac{1}{2} + \frac{i\sqrt{3}}{2}.
\]

Now consider the case $G = S_3$ and $(r,s) = (3,0)$. It turns out that Question~\ref{question:main-q} is already very difficult, but is easier when one allows non-maximal orders. Let $\cF^{\ord}(G,r,s)$ denote the set of (isomorphism classes of) orders in degree $n$ number fields with Galois group $G$ and signature $(r,s)$, and let $\Omega^{\ord}(G,r,s)$ be the corresponding set of points. 

In this case, we have:
\begin{itemize}
    \item work of Cusick \cite{cusick} in 1991 shows that the closure $\overline{\Omega^{\ord}(S_3,3,0)}$ contains an arc;
    \item work of David and Shapira \cite{david-shapira} in 2018 shows that $\overline{\Omega^{\ord}(S_3,3,0)}$ contains an infinite set of curves in $\GL_2(\ZZ) \backslash \mfh$, and conjectured that $\overline{\Omega^{\ord}(S_3,3,0)} = \GL_2(\ZZ)\backslash \mfh$;
    \item work of Dang, Gargava, and Li \cite{dang-gargava-li} in 2025 proves the conjecture of David and Shapira mentioned above;
    \item and independently work of Corso and Hertz \cite{corso-hertz} in 2025 shows that $\overline{\Omega^{\ord}(S_3,3,0)}$ is not compact;
\end{itemize}

To the best of our knowledge, almost nothing is known about \emph{maximal} orders, i.e. $\overline{\Omega(S_3,3,0)}$. The methods in \cite{cusick}, \cite{david-shapira}, \cite{corso-hertz} involve producing orders which are likely maximal, but a very difficult squarefree sieve is needed to \emph{prove} that they are maximal. In contrast, the methods of \cite{dang-gargava-li} involve producing orders which are provably very far from maximal (in fact, they are not even Gorenstein). 
}

Another natural case to consider is when $G$ is equal to the dihedral group $D_p$, for $p$ a prime. Indeed, in \cite{harron-holmes-vemulapalli}, the third author and the last author, joint with Harron, prove that $\overline{\Omega^{\ord}(D_p,1,\frac{p-1}{2})}$ is contained in an explicit finite union of translates of periodic torus orbits. For example, when $p = 5$, the unit lattice of a $D_5$ field with signature $(1,2)$ has rank $2$, and it is shown that the unit shape is contained in the arc 
\[
\left(x+\frac{1}{2}\right)^2+\left(y-\frac{1}{2\sqrt{3}}\right)^2=\left(\frac{2}{\sqrt{3}}\right)^2
\]
from $\frac{1+i\sqrt{3}}{2}$ to $-\frac{1}{2}+i\frac{5}{2\sqrt{3}}$ in $\gl_2(\ZZ)\backslash \mathfrak{h}$.

\bpoint{Main Results}
\label{subsec:results}

In this paper, we investigate Question~\ref{question:main-q} for $D_4$-quartic extensions with $2$ real embeddings and $1$ pair of complex conjugate embeddings. Such fields have a unit lattice of rank $2$, and therefore are a natural starting point for investigations on unit lattices. 

For the rest of this article, let $L \in \cF(D_4,2,1)$. As usual, the space of rank $2$ lattices up to homothety and reflection is identified with $\gl_2(\ZZ)\backslash \mfh$, where here $\mfh$ is the upper half plane. Let $\cF$ be the standard fundamental domain of $\gl_2(\ZZ)\backslash \mfh$, and let $x,y$ be the standard coordinates. Thus, the primary question in this paper is the following. 

\tpoint{Main Question for $D_4$ with signature $(2,1)$} { \em
\label{q:main}
What is $\Omega(D_4,2,1)$?
}

As in the totally real $S_3$ case discussed above, we remark that answering Question~\ref{q:main} for (possibly non-maximal) orders in $D_4$-quartic fields with signature $(2,1)$ is significantly easier than answering the problem for maximal orders.  In this article, we focus on the case of maximal orders and, while the answer to Question~\ref{q:main} remains mysterious, we show that $\Omega(D_4,2,1)$ consists entirely of transcendental elements on the boundary of $\cF$ (see Theorem ~\ref{prop:boundary-and-trans}). Despite the fact that $\Omega(D_4,2,1)$ is entirely transcendental, we further prove that its set of limits points contains at least three explicit algebraic numbers (see Theorem~\ref{thm:limit-points}) and conjecture the following:

\tpoint{Conjecture} { \em
The set of limit points of ${\Omega(D_4,2,1)}$ equals the boundary of $\cF$.
}

Our conjecture is supported by computational evidence; see Figure $1$. It would be interesting to conjecture a distribution for $\{p_L\}_{L \in \cF(D_4,2,1)}$ given some ordering of $\cF(D_4,2,1)$ (three natural choices of ordering are the ordering by absolute discriminant, conductor, and product of ramified primes). 

Our first major result is that the shape of $E_L$ is constrained, and the point $p_L$ associated to it \emph{must} be transcendental. 

\tpoint{Theorem} { \em
\label{prop:boundary-and-trans}
Let $L$ be a $D_4$-quartic field with signature $(2,1)$. Then $p_L \in \cF$ is a transcendental number. Moreover, $p_L$ is contained in the boundary of $\cF$, i.e. is contained in either
\begin{enumerate}
    \item the left vertical line $x = 0$ and $1 \leq y$,
    \item the lower arc $1 = x^2 + y^2$ for $0 \leq x \leq 1/2$ and $y > 0$,
    \item or the right vertical line $x = 1/2$ and $\frac{\sqrt{3}}{2} \leq y$.
\end{enumerate}
}

In particular, $p_L$ is never equal to the cuspidal points $i$ and $\frac{1}{2} + \frac{i\sqrt{3}}{2}$. One basic question one may ask about the topology of $\Omega(D_4,2,1)$ is: what is its set of limit points?\footnote{In a topological space $X$, we say $x$ is a limit point of a set $S \subseteq X$ if for every open set $x \in U$, there exists a point $y \in S \cap U$ with $y \neq x$.}. In fact, it is not too difficult to show that there are infinite subsets of $\Omega(D_4,2,1)$ with no limit points. 

\tpoint{Theorem} { \em
\label{thm:escape-of-mass}
Given a real quadratic field $K$, define $\cF(K)$ the fields in $\cF(D_4,2,1)$ which contain $K$ as a subfield. Then the multiset $\{p_L\}_{L \in \cF(K)}$ has no limit points. \footnote{A limit point of the multiset $\{p_L\}_{L\in \cF(K)}$ is a point $p$ and an infinite sequence $\{L_i\}$ of distinct fields such that $p_{L_i} \rightarrow p$.}
}

\tpoint{Corollary} {\em
\label{cor:escape-of-mass}
With the notation as in Theorem~\ref{thm:escape-of-mass}, there are only finitely many points of the multiset $\{p_L\}_{L \in \cF(K)}$ which lie on the lower arc on the boundary of $\cF$.
}

To prove Corollary~\ref{cor:escape-of-mass} from Theorem~\ref{thm:escape-of-mass}, assume for the sake of contradiction there are infinitely many fields $L$ such that $p_L$ lies on the lower arc. If there exists a point $p$ such that $p = p_L$ for infinitely many $L$, then $p$ is a limit point of $\{p_L\}_{L \in \cF(K)}$, which is a contradiction. Therefore, the set of points $p$ which lie on the lower arc and are equal to $p_L$ for some $L$ is infinite. Because the lower arc is compact, there exists some convergent sequence $\{p_L\}$ and hence a limit point, contradiction.

However, when the real quadratic subfield $K$ is allowed to vary, it turns out that $\Omega(D_4,2,1)$ \emph{does} contain limit points, and these limit points may even be algebraic. Namely:

\tpoint{Theorem} { \em
\label{thm:limit-points}
The set of limit points of $\Omega(D_4,2,1)$ contains $\left\{i\sqrt{3}, \frac{1}{2} + \frac{i\sqrt{3}}{2}, \frac{1}{7} + \frac{4i\sqrt{3}}{7}\right\}$.
}

See Figure~\ref{fig:fig-5} for a picture of these points. Note that the last point, $\frac{1}{7} + \frac{4i\sqrt{3}}{7}$, lies on the interior of the lower arc of the boundary of $\cF$. Thus we obtain the following corollary.

\tpoint{Corollary} {\em
Infinitely many points of $\Omega(D_4,2,1)$ lie on the lower arc in the boundary of $\cF$.
}

\bpoint{Motivation from cryptography}
\label{subsection:motivation}
Over the past decades, lattices have gained increasing attention due to their significance in post-quantum cryptography. The hardness of certain lattice problems plays a critical role in ensuring the security of cryptographic schemes against both classical and quantum attacks. One of the most promising candidates for lattice-based cryptography is the Learning With Errors (LWE) problem \cite{Regev}, along with its variants. At a high level, LWE involves recovering a secret from a system of linear equations perturbed by some noise, drawn from a specific error distribution. The strong security guarantees of LWE-based schemes rely on reductions to worst-case lattice problems such as the Shortest Vector Problem (SVP) and the Closest Vector Problem (CVP).

The shape of a lattice describes its geometric structure, which influences the
complexity of these fundamental problems. In addition, the security of LWE–based
schemes depends on the choice of the error distribution. Any leakage of
information about this distribution could be exploited by attackers \cite{ELOS}. Therefore, understanding the shape of lattices can inform our choices of secure parameters and may help identify weak instances in lattice-based cryptographic schemes.

Extending our analysis to higher-rank cases is an interesting direction for future research as well. 

\bpoint{Acknowledgments}
\label{subsec:acknowledgements}
We are grateful to Heidi Goodson, Allechar Serrano L\'opez, and Mackenzie West for organizing the Rethinking Number Theory conference, where this project began. We also thank Ha Tran for helpful comments. We also would like to thank Jose Miguel Cruz Rangel for help understanding Stender's theorem.  We are also grateful to the anonymous referee for helpful comments. S. Vemulapalli was funded by the NSF under grant number DMS2303211.

\section{Notation}
\label{sec:notation}
As before, let $\cF(D_4,2,1)$ be the set of isomorphism classes of $D_4$-quartic extensions with two real embeddings and one pair of complex conjugate embeddings. Let $L$ denote an element of $\cF(D_4,2,1)$ and let $K$ be the (unique) real quadratic subfield contained in $L$. Let $\sigma \in \Gal(L/K)$ be the nontrivial automorphism of $L$ over $K$. As usual, let $E_L$ denote the unit lattice of $L$ and let $\iota \colon \cO_L^{\times} \rightarrow E_L$ denote the quotient. Let $\Log \colon E_L \rightarrow \RR^3$ denote the logarithmic Minkowski embedding, and let $\sigma_1,\sigma_2$ be real embeddings of $L$, let $\tau$ be a choice of complex conjugate embedding of $L$, and let $\overline{\tau}$ be its complex conjugate. For any $x \in E_L$, let $\lvert x \rvert$ denote the length of $x$, with respect to the Euclidean metric, in the embedding $\Log \colon E_L \rightarrow \RR^{3}$.

The element $\sigma \in \Gal(L/K)$ acts on $\cO_L^{\times}$ and the action descends to the quotient $E_L = \cO_L^{\times}/\{\pm 1\}$. By abuse of notation, we think of $E_L$ as an additive group isomorphic to $\ZZ^2$ equipped with an action of $\sigma$. Therefore, for any two units $u,v \in \cO_L^{\times}$, we write $\iota(uv) = \iota(u) + \iota(v)$.

As before, the space of rank $2$ lattices up to homothety and reflection is identified with $\gl_2(\ZZ)\backslash \mfh$, where here $\mfh$ is the upper half plane. Let $\cF$ be the standard fundamental domain of $\gl_2(\ZZ)\backslash \mfh$, and let $x,y$ be the standard coordinates.

\subsection{Vinogradov notation}
Throughout the paper we use the following asyptotic notation, first introduced by Vinogradov. For a subset $S \subseteq \RR_{\geq 0}$ and two functions $f, g \colon S \rightarrow \RR$ and some parameter $\epsilon$, the notation $f(x) \ll_{\epsilon} g(x)$ means that there exist positive real numbers $M_{\epsilon}$ and $C_{\epsilon}$, dependent on $\epsilon$, such that for all $x \in S$ with $x \geq M_{\epsilon}$ we have $f(x) \leq C_{\epsilon}g(x)$. Similarly, the notation $f(x) \asymp_{\epsilon} g(x)$ means that there exist positive real numbers $M_{\epsilon}$, $C_{\epsilon}$, and $D_{\epsilon}$, dependent on $\epsilon$, such that for all $x \in S$ with $x \geq M_{\epsilon}$ we have $C_{\epsilon}g(x) \leq f(x) \leq D_{\epsilon}g(x)$.

We also abuse the big-O and little-o notation in the following way. We write $f(x) = g(x) + O(h(x))$ if there exist positive real numbers $M,C$ such that for all $x \in S$ with $x \geq M$, we have $\lvert f(x) - g(x) \rvert \leq C h(x)$. Similarly, we write $f(x) = g(x) + o(h(x))$ if for every positive real number $C$, there exists a positive real number $M_C$ such that for all $x \in S$ with $x \geq M_C$, we have $\lvert f(x) - g(x)\rvert \leq Ch(x)$.

\section{Proof of Proposition~\ref{prop:boundary-and-trans}}
\label{sec:proof-boundary-and-trans}

	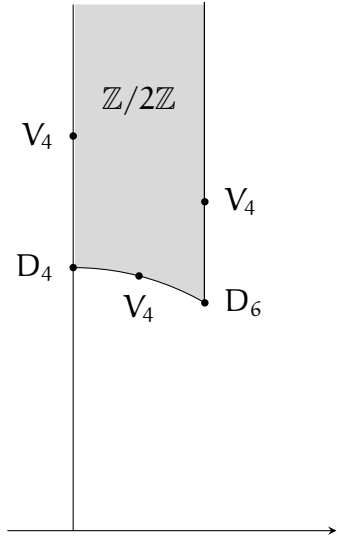
\begin{figure}
        \label{fig:fig-5}
	\centering
		\begin{tikzpicture}[scale = 3.5]
			\tikzset{myptr/.style={decoration={markings,mark=at position 1 with %
    							{\arrow[scale=3,>=stealth]{>}}},postaction={decorate}}}
		
			\begin{scope}
    				\clip (0.007,.86) rectangle (.5,2.01);
    				\shade[
					bottom color = lightgray!60,
					middle color = lightgray!60,
					top color = lightgray!60,
					shading angle = 0
					] (0,1) -- (.5,.866) -- (.5,2) -- (0,2);
					
				\draw[thick] (.5,.866) -- (.5,2.01);
				\filldraw[fill=white, draw=black] (0,0) circle (1);
				
			\end{scope} 
            \draw[-] (0,0)--(0,2) ;
			\draw[->, >=stealth] (-.25,0)--(1,0);
            
			\node[scale=.75, label = $\ZZ/2\ZZ$] (I) at (.25, 1.5) {};
            \node[scale=.75, label = left:$D_4$] (LC) at (0, 1) {};
            \node[scale=.75, label = right:$D_6$] (RC) at (.5, .866) {};
            \node[scale=.75, label = left:$V_4$] (LB) at (0, 1.5) {};
            \node[scale=.75, label = right:$V_4$] (RB) at (.5, 1.25) {};
            \node[scale=.75, label = below:$V_4$] (BB) at (.25, .968) {};

            \node at (LB) [circle,fill,inner sep=1pt]{};
            \node at (RB) [circle,fill,inner sep=1pt]{};
            \node at (LC) [circle,fill,inner sep=1pt]{};
            \node at (RC) [circle,fill,inner sep=1pt]{};
            \node at (BB) [circle,fill,inner sep=1pt]{};
            
		\end{tikzpicture}
		\caption{Each point of $\cF$ corresponds to an equivalence class of rank $2$ lattices up to homothety and reflection. The groups are the automorphism groups of the corresponding lattices; for example, the cuspidal point $i$ corresponds to a lattice with automorphism group $D_4$. A lattice on the interior of $\cF$ has only one nontrivial automorphism, given by multiplication by $-1$.}
		\label{rank2shapes}
	\end{figure}

Proposition~\ref{prop:boundary-and-trans} will follow from Proposition~\ref{prop:boundary} and Proposition~\ref{prop:trans}.

\tpoint{Proposition} { \em
\label{prop:boundary}
For any $L \in \cF(D_4,2,1)$, the point $p_L$ is contained in the boundary of $\cF$. 
}

\bpf
Together, Lemma~\ref{lem:nontriv-aut} and Lemma~\ref{lem:preserve-length} imply that $E_L$ is a rank $2$ lattice with an order two automorphism which is not given by multiplication by $-1$. It follows from standard results that the only such rank $2$ lattices, up to homothety, correspond to the boundary of $\cF$.
\epf

Recall that by abuse of notation, we think of $E_L$ as an additive group isomorphic to $\ZZ^2$ equipped with an action of $\sigma$; see Section~\ref{sec:notation}.

\tpoint{Lemma} {\em
\label{lem:nontriv-aut}
The action of $\sigma$ on $E_L$ is not given by multiplication by $\pm 1$.
}

\bpf
The rank $2$ lattice $E_L$ has a rank $1$ sublattice $E_K$; clearly $\sigma$ fixes $E_K$, so the action of $\sigma$ on $E_L$ is not given by multiplication by $-1$. 

Suppose for the sake of contradiction that $\sigma$ acts trivially on $E_L$. Choose $u \in \cO_L^{\times}$ such that $n\iota(u) \notin E_K$ for any nonzero integer $n$. By assumption $\sigma(\iota(u)) = \iota(u)$. Because the kernel of $\iota$ is given by roots of unity, this implies that
\begin{equation}
    \frac{\sigma(u)}{u} = \zeta
\end{equation}
for some root of unity $\zeta$. However, the only roots of unity in $L$ are $\pm 1$, so either $\sigma(u) = u$ or $\sigma(u) = -u$. 

Now, $\sigma(u) \neq u$ because $u \notin K$. So $\sigma(u) = -u$. Since $L/K$ is a quadratic extension, we may write $L=K(\sqrt{\alpha})$ for some $\alpha \in K$ and $u=a+b\sqrt{\alpha}$ for $a,b \in K$. Then $\sigma(u) = a - b\sqrt{\alpha}$. Because $\sigma(u) = -u$, we have that $u = b\sqrt{\alpha}$; but then $u^2 = b^2 \alpha \in K$ which implies that $2\iota(u) \in E_K$, which is a contradiction. 
\epf

\tpoint{Lemma} {\em
\label{lem:preserve-length}
For any $x \in E_L$, we have that $\lvert \sigma(x) \rvert = \lvert x \rvert$.
}

\bpf
We have $\sigma_1(\sigma(x)) = \sigma_2(x)$ and $\sigma_2(\sigma(x)) = \sigma_1(x)$ and $\tau(\sigma(x)) = \overline{\tau(x)}$. Therefore:
\[
    \lvert \sigma(x) \rvert =  \sqrt{\lvert \sigma_1(\sigma(x)) \rvert^2 + \lvert \sigma_1(\sigma(x)) \rvert^2 + 4\lvert \tau(\sigma(x)) \rvert^2} = \sqrt{\lvert \sigma_2(x) \rvert^2 + \lvert \sigma_1(x) \rvert^2 + 4\lvert \overline{\tau(x)} \rvert^2} = \lvert x \rvert.
\]
\epf

We now prove the second half of Theorem~\ref{prop:boundary-and-trans}, i.e. that $p_L$ is transcendental.

\tpoint{Proposition} { \em
\label{prop:trans}
For any $L \in \cF(D_4,2,1)$, the point $p_L$ is transcendental.
}

\bpf 
Suppose for the sake of contradiction that $p_L$ is algebraic.  Let $\cL= \ZZ \langle 1,p_L \rangle$ be the sublattice of $\CC$ spanned by $1$ and $p_L$. Because $p_L$ is algebraic, the length of any vector in $\cL$ is algebraic, so the quotient of the lengths of any two nonzero vectors in the lattice $\cL$ is also algebraic. 

Now, $E_L$ is homothetic to $\cL$, so this implies that the quotient of the lengths of any two nonzero vectors in $E_L$ is algebraic. We now will exhibit two distinct vectors in $E_L$ whose length quotient is not algebraic, thus showing a contradiction and proving that $p_L$ is transcendental. 

\noindent \textbf{Step 1: defining $u_1$ and $u_2$.}
Choose $u_1 \in \cO_L^{\times}$ so that $\iota(u_1) \neq 0$ and $\iota(u_1) \in \ker(N_{L/K} \colon E_L \rightarrow E_K)$. 

Now, because $\iota(u_1) \in \ker(N_{L/K} \colon E_L \rightarrow E_K)$, we have that $\iota(u_1) + \iota(\sigma((u_1)) = 0$, so
\begin{equation}
\label{eqn:norm}
    \Log(\iota(u_1)) + \Log(\iota(\sigma(u_1)))  = (0,0,0).
\end{equation}
Write:
\begin{equation}
\label{eqn:u1}
     \Log(\iota(u_1)) = (\log \lvert \sigma_1(u_1) \rvert, \log \lvert \sigma_2(u_1) \rvert, 2\log \lvert \tau(u_1) \rvert)
\end{equation}
Recall the action of $\sigma$ on the embeddings; $\sigma_1 \circ \sigma = \sigma_2$, $\sigma_2 \circ \sigma = \sigma_1$ and $\tau \circ \sigma = \overline{\tau}$. So:
\begin{equation}
\label{eqn:sigma-u1}
\Log(\iota(\sigma(u_1))) = (\log \lvert \sigma_2(u_1) \rvert, \log \lvert \sigma_1(u_1) \rvert, 2\log \lvert \tau(u_1) \rvert)
\end{equation}

Substituting Equation~\ref{eqn:u1} and Equation~\ref{eqn:sigma-u1} into Equation~\ref{eqn:norm}, we obtain that $\Log \lvert \tau(u_1) \rvert = 0$. Finally, the image of $\Log$ lies on the hyperplane whose coordinates sum to zero, so we may rewrite Equation~\ref{eqn:u1} as
\begin{equation}
\label{eqn:u1-new}
     \Log(\iota(u_1)) = (\log \lvert \sigma_1(u_1) \rvert, -\log \lvert \sigma_1(u_1) \rvert, 0).
\end{equation}

So, 
\begin{equation}
\label{eqn:u1-length}
\lvert \iota(u_1) \rvert = \sqrt{2 (\log \lvert \sigma_1(u_1) \rvert)^2} = \Big \lvert \sqrt{2} \log \lvert \sigma_1(u_1) \rvert \Big \rvert
\end{equation}

Now choose $u_2 \in \cO_K^{\times}$ so that $\iota(u_2) \neq 0$. Write:
\begin{equation}
\label{eqn:u2}
     \Log(\iota(u_2)) = (\log \lvert \sigma_1(u_2) \rvert, \log \lvert \sigma_2(u_2) \rvert, 2\log \lvert \tau(u_2) \rvert)
\end{equation}
Because $\sigma$ fixes $u_2$, we have that $\log \lvert \sigma_1(u_2) \rvert = \log \lvert \sigma_2(u_2) \rvert$. Because the image of $\Log$ lies on the hyperplane whose coordinates sum to zero, we obtain
\begin{equation}
\label{eqn:u2-new}
     \Log(\iota(u_2)) = (\log \lvert \sigma_1(u_2) \rvert, \log \lvert \sigma_1(u_2) \rvert, -2\log \lvert \sigma_1(u_2) \rvert).
\end{equation}

\begin{equation}
\label{eqn:u2-length}
\lvert \iota(u_2) \rvert = \sqrt{6 (\log \lvert \sigma_1(u_2) \rvert)^2} = \Big \lvert \sqrt{6} \log \lvert \sigma_1(u_2) \rvert \Big \rvert
\end{equation}
Combining Equation~\ref{eqn:u1-length} and Equation~\ref{eqn:u2-length}, we get that
\begin{equation}
\label{eqn:ratio}
\frac{\lvert \iota(u_1) \rvert}{\rvert \iota(u_2) \rvert} = \frac{\big \lvert \log \lvert \sigma_1(u_1) \rvert \big \rvert}{\sqrt{3}\big \lvert \log \lvert \sigma_1(u_2) \rvert\big \rvert}
\end{equation}

\noindent \textbf{Step 2: showing $\log \lvert \sigma_1(u_1) \rvert $ and $ \log \lvert \sigma_1(u_2) \rvert$ are $\QQ$-linearly independent complex numbers. }
Suppose for the sake of contradiction there existed integers $a,b$ such that:
\[
    a \log \lvert \sigma_1(u_1) \rvert = b \log \lvert \sigma_1(u_2) \rvert
\]
then
\[
    \log \lvert \sigma_1(u_1^a) \rvert = \log \lvert \sigma_1(u_2^b) \rvert.
\]
Because $\sigma_1$ is a real embedding, this implies that
\[
    u_1^a = \pm u_2^b.
\]
Because $u_2 \in \cO_K^{\times}$ this implies that $u_1^a \in \cO_K^{\times}$. Then because $\sigma$ fixes $K$, we have that
\begin{equation}
\label{eqn:u1-in-K}
\log \lvert \sigma_1(u_1) \rvert = \log \lvert \sigma_2(u_1) \rvert.
\end{equation}
Combining Equation~\ref{eqn:u1-in-K} with Equation~\ref{eqn:u1-new} implies that $\Log(u_1) = 0$; but this implies $\iota(u_1) = 0$, a contradiction.

\noindent \textbf{Step 3: Applying linear independence of logarithms.} Because $\log \lvert \sigma_1(u_1) \rvert $ and $ \log \lvert \sigma_1(u_2) \rvert$ are $\QQ$-linearly independent, the Gelfond-Schneider theorem on linear independence of logarithms (\cite{gelfand, schneider}) implies that they are $\overline{\QQ}$-linearly independent. Combining this  Equation~\ref{eqn:ratio}, we obtain that $\lvert \iota(u_1) \rvert /\lvert \iota(u_2) \rvert$ is transcendental.
\epf

\section{Proof of Theorem~\ref{thm:escape-of-mass}}
\label{sec:proof-escape-of-mass}

The key input is the following lower bound on the regulator of a number field due to Silverman, in terms of its discriminant. We remark that this bound has recently been improved by Akhtari and Vaaler (\cite{akhtari-vaaler}). For a number field $L$, let $\Reg(L)$ denote its regulator. 

\tpoint{Theorem 1, \cite{silverman}} { \em
For any number field $L$, we have
\[
	\Reg(L) > 2^{-4d^2} \big(\log(d^{-d^{\log_2 8d}}\lvert\Disc(L)\rvert)\big)^{r(L) - \rho(L)}
\]
where $r(L)$ is the rank of the (free part of the) unit group of $L$ and 
\[
    \rho(L) = \max\{r(K) : K \subset L, \; K \neq L\}.
\]
}

\tpoint{Proposition} { \em
\label{prop:y-bound}
Let $L$ be a $D_4$-quartic field with signature $(2,1)$ containing a real quadratic field $K$. Write $p_L = x_L + iy_L$. Then there exists a positive real number $c_K$ dependent only on $K$ such that
\[
	y_L \geq c_K \sqrt{\log \lvert \Nm_{K/\QQ}(\Disc(L/K)) \rvert}.
\]
}

\bpf
Let $u,v \in E_L$ be the two shortest linearly independent vectors of $E_L$, ordered so that $\lvert u \rvert \leq \lvert v \rvert$; then $\{u,v\}$ is a reduced basis of $E_L$. Let $\epsilon \in E_K$ be a generator. Let $\theta$ be the angle between $u$ and $v$. Because $\{u,v\}$ is a reduced basis of $E_L$, we have $\lvert \sin \theta \rvert \geq \frac{\sqrt{3}}{2}$.
Then
\[
    y = \frac{\lvert v \rvert}{\lvert u \rvert}\lvert\sin \theta \rvert \geq \frac{\sqrt{3}}{2} \frac{\lvert v \rvert}{\lvert u \rvert} \geq \frac{\sqrt{3}}{2} \frac{\lvert v \rvert}{\lvert \epsilon \rvert} = \frac{\sqrt{3}}{2\sqrt{2}} \frac{\lvert v \rvert}{\Reg(K)}.
\]
Now, we have that
\[
    \lvert v \rvert^2 \geq \lvert u \rvert \lvert v \rvert \geq \Covol(E_L) = \sqrt{3} \Reg(L),
\]
so
\[
    \lvert v \rvert \geq 3^{1/4}\sqrt{\Reg(L)}.
\]
Combining, we get that
\[
    y_L \geq \frac{3^{3/4}}{2^{3/2}}\frac{\sqrt{\Reg(L)}}{\Reg(K)}  \gg \frac{\sqrt{\log(\lvert\Disc(L)\rvert)}}{\Reg(K)} \gg_K \sqrt{\log(\lvert\Disc(L)\rvert)}.
\]
The discriminant of $L$ factors as:
\[
    \Disc(L) = \Disc(K)^2 \Nm_{K/\QQ}(\Disc(L/K)),
\]
and thus we obtain
\[
    y_L \gg_{K} \sqrt{\log \lvert \Nm_{K/\QQ}(\Disc(L/K)) \rvert}.
\]
\epf

Now, Theorem~\ref{thm:escape-of-mass} follows immediately from the fact that for every positive real number $X$, there are finitely many quadratic extensions $L$ of $K$ such that $\lvert \Nm_{K/\QQ}(\Disc(L/K)) \rvert \leq X$.

\section{Proof of Theorem~\ref{thm:limit-points}}
\label{sec:proof-limit-points}

Theorem~\ref{thm:limit-points} follows from Proposition~\ref{prop:first-point}, Proposition~\ref{prop:second-point}, and Proposition~\ref{prop:third-point}. For reference, we include a table of the explicit families used in this section, Table \ref{tab:family_table}, as well as a figure illustrating the three limit points, Figure \ref{limitpoints}:

\begin{table}[htbp]
\centering
\begin{tabular}{|c|c|c|c|}
\hline
\textbf{Citation} 			&	 \textbf{Family:} $\{f_n(x)\}$ 		       & 	\textbf{Units: $\{u_1, u_2\}$} 				                             & 	\textbf{Limit Point} 			            \\
\hline
                    		& 	                           	               & 	                                                                         &    	 				                             \\
Stender \cite{stender2} 		& 	$\{x^4 - \omega\}$ 	                       & 	$\left\{\frac{\omega+n}{\omega-n}, \frac{n^3}{(\omega-n)^4}\right\}$    	 & 	$i\sqrt{3}$ 				                 \\ 
                    		& 	 \small{where $\omega =(n^4+n^3)$}	       & 	                                                                         & 	 				                             \\
 \hline
                     		& 	                          	               & 	                                                                             & 	 				                             \\
Stender \cite{stender2} 		& 	$\{x^4-\omega\}$ 		                   & 	$\left\{\frac{\omega+n}{\omega-n}, \frac{1}{\omega-n}\right\}$ 		          & 	$\frac{1}{7}+\frac{4i\sqrt{3}}{7}$            \\ 
                    		& 	 \small{where $\omega =(n^4-1)$}	       & 	                                                                              & 	 				                                 \\
\hline
                    		& 	                           	               & 	     $\{\rho, \epsilon\}$                                                     & 	 				                             \\
LPS \cite{lps} 	        & 	$\{x^4-nx^3+(1-n)x^2-2x+1\}$ 	           & 	    \small{$\rho$ a root of $f_n(x)$}                                                 & 	$\frac{1+i\sqrt{3}}{2}$              \\
                    		& 	                                           & 	\small{$\epsilon = \frac{-n+\sqrt{n^2-4}}{2}$}                                      & 	 				                              \\
\hline
                    		& 	                           	               & 	                                                                              & 	 				                             \\
Nakamula \cite{nk} 			& 	$\{x^4-nx^3+3x^2-nx+1\}$ 		           & 	$\left\{\rho, \rho+\rho^{-1}\right\}$ 			                              & 	$i\sqrt{3}$                                   \\
                    		& 	                                           & 	\small{$\rho$ is a root of $f_n(x)$}                                                  & 	 				                                 \\
\hline
\end{tabular}
\caption{Summary of polynomial families, unit group generators, and limit points for $D_4$-quartic fields with signature $(2,1)$.}
\label{tab:family_table}
\end{table}

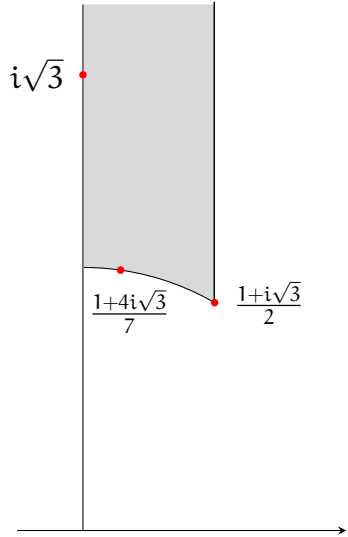
\begin{figure}[htbp]
\centering
    \begin{tikzpicture}[scale = 3.5]
        \tikzset{myptr/.style={decoration={markings,mark=at position 1 with %
                            {\arrow[scale=3,>=stealth]{>}}},postaction={decorate}}}
    
        \begin{scope}
                \clip (0.007,.86) rectangle (.5,2.01);
                \shade[
                bottom color = lightgray!60,
                middle color = lightgray!60,
                top color = lightgray!60,
                shading angle = 0
                ] (0,1) -- (.5,.866) -- (.5,2) -- (0,2);
                
            \draw[thick] (.5,.866) -- (.5,2.01);
            \filldraw[fill=white, draw=black] (0,0) circle (1);
            
        \end{scope} 
        \draw[-] (0,0)--(0,2) ;
        \draw[->, >=stealth] (-.25,0)--(1,0);

        \node[scale=.75, label = left:$i\sqrt{3}$]                       (stender1) at (0, 1.732) {};
        \node[scale=.75, label = below:$\;\;\frac{1+4i\sqrt{3}}{7}$]          (stender2) at (1/7, 4*1.732/7) {};
        \node[scale=.75, label = right:$\frac{1+i\sqrt{3}}{2}$]     (lps) at (.5, .866) {};

        \node at (stender1) [circle,fill,color=red,inner sep=1pt]{};
        \node at (stender2) [circle,fill,color=red,inner sep=1pt]{};
        \node at (lps)      [circle,fill,color=red,inner sep=1pt]{};

    \end{tikzpicture}
    \caption{Limit points of the unit lattices of the families of quartic fields arising from constructions of Stender, Lepr\'evost-Pohst-Sch\"opp, and Nakamula. }
    \label{limitpoints}
\end{figure}

To show that there exists a sequence of $D_4$-quartic fields whose unit lattices have shape converging to the points $i\sqrt{3}$ and $\frac{1 + 4i\sqrt{3}}{7}$, we make use of an explicit construction of Stender \cite{stender2} describing the units of some pure quartic fields. Let $n$ and $d$ be positive integers such that $n^4 \pm d$ is squarefree, $d$ is not a square, $d \mid n^3$, and $\frac{n^4}{d} \pm 1$ is squarefree. Write $\omega = (n^4 \pm d)^{1/4}$ and set $L = \QQ(\omega)$. It is easy to see that $L$ is a $D_4$-quartic field with signature $(2,1)$. 

\tpoint{Theorem (Stender, \cite{stender2})} { \em
\label{thm:stender}
The units
\[
   \mathbf{u}= \iota\bigg( \frac{\omega+n}{\omega-n}\bigg) \quad \text { and } \quad \mathbf{v}= \begin{cases}\iota\big(\frac{d}{(\omega-n)^4}\big) & \text { if } d \text { is not a square, } \\
   \iota\big(\frac{\sqrt{d}}{(\omega-n)^2}\big) & \text { if } d \text { is a square and } d \neq 1, \\  \iota \big (\frac{1}{\omega-n}\big) & \text { if } d=1\end{cases}
\]
form a basis of $E_{L}$, except for the three cases $\omega^4=8, \omega^4=12$ and $\omega^4=20$.}
\bpoint{The point $i\sqrt{3}\in \Omega(D_4, 2,1)$}
\tpoint{Proposition} { \em
\label{prop:first-point}
The set of limit points of $\Omega(D_4,2,1)$ contains $i\sqrt{3}$.
}

As Stender's original theorem is in German, we now give a description of the key ideas behind this proof. (Where appropriate, we've modernized the proof to simplify exposition). 

Given a $D_4$-quartic extension $L/\QQ$ with quadratic subfield $K$ and one real embedding, there is a natural norm map
\[
    \Nm : E_L \rightarrow E_K.
\]
The kernel has rank $1$; let $\epsilon_0$ be a generator of the kernel (chosen appropriately). Let $\epsilon_K$ be a generator of $E_K$. Depending on whether $p_L$ lies on the left boundary, the bottom boundary, or the right boundary, there are three possibilities for a basis of $E_L$. To see this, identify $E_L$ with the lattice in $\RR^2$ spanned by $p_L$ and $1$ and let $u,v$ be units corresponding to $p_L$ and $1$ respectively. In all three cases, we write down an explicit formula for the norm map on $E_L$, identified with $\ZZ\langle1, p_L\rangle$ and use that norm map to explicitly write $\epsilon_0$ and $\epsilon_K$ as a function of $u$ and $v$. This explicit formula allows us to produce a basis of $E_L$, as as function of $\epsilon_0$ and $\epsilon_K$. Namely:
\begin{itemize}
    \item Suppose $p_{L} = iy_0$ lies on the left boundary. Then the norm map is reflection across the $x$-axis, so it sends $x + iy \rightarrow 2x$. So, in this case $\epsilon_0 = u$, and $\epsilon_K = v$, so $\{\epsilon_0,\epsilon_K\}$ is a basis of $E_L$;
    \item if $p_L$ lies on the bottom boundary, then the norm map reflects along the line bisecting $p_L$ and $1$. Concretely, $\epsilon_0 = uv^{-1}$ and $\epsilon_K = uv$. Therefore, $\{\sqrt{\epsilon_0\epsilon_K},\sqrt{\epsilon_0^{-1}\epsilon_K}\}$ is a basis of $E_L$.
    \item if $p_L$ lies on the right boundary, then the norm map is again reflection across the $x$-axis, so $\epsilon_0 = u^2v^{-1}$ and $\epsilon_K = v$. Therefore, $\{\epsilon_K, \sqrt{\epsilon_0\epsilon_K}\}$ is a basis of $E_L$.
\end{itemize}

The rest of the proof of Stender's theorem has two steps. First, Stender determines which of the three boundary components $p_L$ lies in from $\omega$, $n$, and $d$. Next, Stender explicitly calculates $\epsilon_0$ and $\epsilon_K$ in terms of $\omega$, $n$, and $d$. (In fact, calculating $\epsilon_K$ was done earlier by G. Degert \cite{Degert1958}). In order to obtain $\epsilon_0$, Stender proceeds as follows. It was previously proved (by Bernstein and Hasse) that $\zeta = \frac{\omega + n}{\omega - n}$ is a unit. To obtain $\epsilon_0$, Stender uses the Bernstein--Hasse result that $\zeta = \frac{\omega+n}{\omega-n}$ is a unit and shows that $\zeta$ is the generator of the kernel of the norm map $E_L \to E_K$.

Stender proves that in fact this special unit is a generator for the kernel of the norm map $E_L \rightarrow E_K$. 

\epoint{Proof of Proposition~\ref{prop:first-point}}

\noindent \textbf{Step 1: Rewriting the basis.}

We now compute the basis in the special case where $d=n^3$ and $n$ is an odd prime number. Write $\omega=(n^4+d)^{1/4}=(n^4+n^3)^{1/4}$ and set ${L_n}= \mathbb{Q}(\omega)$. Note that there are infinitely many primes $n$ such that $n + 1$ is squarefree; see Theorem 2 of \cite{mirsky}. By Theorem~\ref{thm:stender} 
the units
$
    \left\{\iota\left(\frac{\omega+n}{\omega-n}\right), \iota\left(\frac{n^3}{(\omega-n)^{4}}\right)\right\}
$
form a basis of $E_{L_n}$.

Let $\mathbf{u} = \frac{\omega+n}{\omega-n}$ and $\mathbf{v} = \frac{n^3}{(\omega - n)^4}$. We compute the logarithmic embeddings.
\begin{equation}
\label{eqn:log-u}
\Log(\iota(\mathbf{u})) = (\log\lvert \omega + n\rvert-\log \lvert \omega - n\rvert, \log\lvert \omega - n\rvert-\log\lvert\omega + n\rvert,0)
\end{equation}
\begin{equation}
\label{eqn:log-v}
\Log(\iota(\mathbf{v})) = (\log\lvert n^3\rvert-4\log\lvert \omega - n\rvert, \log\lvert n^3 \rvert-4\log\lvert \omega + n \rvert, 2\log\lvert n^3 \rvert -8\log\lvert i\omega + n \rvert)
\end{equation}

A computation shows that:
\begin{equation}
\label{eqn:dot-product}
\Log(\iota(\mathbf{u})) \cdot \Log(\iota(\mathbf{v})) = 4(\log \lvert w + n\rvert - \log \lvert w - n\rvert)^2 = 2\Log(\iota(\mathbf{u})) \cdot \Log(\iota(\mathbf{u}))
\end{equation}

Let $\mathbf{v}_1 = \mathbf{v}\mathbf{u}^{-2}$. Equation~\ref{eqn:dot-product} shows that $\Log(\iota(\mathbf{v}_1)) \cdot \Log(\iota(\mathbf{u})) = 0$, so $\{\iota(\mathbf{u}),\iota(\mathbf{v}_1)\}$ is an orthogonal basis of $E_L$. Therefore, if $\lvert \iota(\mathbf{u}) \rvert \leq \lvert \iota(\mathbf{v}_1) \rvert$, then
\begin{equation}
\label{eqn:pl}
p_L = i\frac{\lvert \iota(\mathbf{v}_1) \rvert}{\lvert \iota(\mathbf{u}) \rvert}.
\end{equation}
 
\noindent \textbf{Step 2: Computing $\lvert \iota(\mathbf{u}) \rvert$ and $\lvert \iota(\mathbf{v}_1) \rvert$.} 

We now compute the numerator and denominator in Equation~\ref{eqn:pl}.. In this case, we have:
\begin{equation}
\label{eqn:taylor-series}
\omega - n = (n^4 + n^{3})^{1/4} - n
= n((1 + n^{-1})^{1/4} - 1) 
\asymp n \cdot n^{-1} 
= 1.
\end{equation}
It is clear that:
\begin{equation}
\label{eqn:bound-1}
\log\lvert \omega + n \rvert = \log n + O(1)
\end{equation}
\begin{equation}
\label{eqn:bound-2}
\log\lvert i\omega + n \rvert = \log n + O(1)
\end{equation}
\begin{equation}
\label{eqn:bound-3}
\log\lvert n^3 \rvert = 3\log n.
\end{equation}
Using Equation~\ref{eqn:taylor-series}, we see that:
\begin{equation}
\label{eqn:bound-4}
\Big \lvert \log\lvert \omega - n \rvert \Big \rvert = O(1). 
\end{equation}

Now, substituting Equation~\ref{eqn:bound-1}, Equation~\ref{eqn:bound-2}, Equation~\ref{eqn:bound-3}, and Equation~\ref{eqn:bound-4} into Equation~\ref{eqn:log-u}, we obtain:
\begin{equation}
\label{eqn:length-u}
  \lvert \iota(\mathbf{u}) \rvert^2 = 2(\log n)^2 + O(\log n).
\end{equation}

Now, $\Log(\iota(\mathbf{v}_1))$ is:
\begin{equation}
\label{eqn:log-v1}
 (\log \lvert n^3 \rvert-2\log \lvert \omega - n\rvert - 2\log\lvert\omega + n \rvert, \log\lvert n^3 \rvert-2\log\lvert \omega + n \rvert - 2\log \lvert \omega - n \rvert, 2\log\lvert n^3 \rvert -8\log\lvert i\omega + n \rvert).
\end{equation}

Substituting Equation~\ref{eqn:bound-1}, Equation~\ref{eqn:bound-2}, and Equation~\ref{eqn:bound-3} into Equation~\ref{eqn:log-v1}, we get that:
\begin{equation}
\label{eqn:log-v1-2}
    \Log(\iota(\mathbf{v}_1)) = (\log n + O(1), \log n + O(1), -2\log n + O(1))
\end{equation}

Squaring the components of Equation~\ref{eqn:log-v1-2}, we get
\begin{equation}
\label{eqn:length-v1}
 \lvert \iota(\mathbf{v}_1) \rvert ^2 = 6(\log n)^2 + O(\log n).
\end{equation}

Substituting Equation~\ref{eqn:length-u} and Equation~\ref{eqn:length-v1} into Equation~\ref{eqn:pl} shows that
\[
	\lim_{p \rightarrow \infty}p_L = i\lim_{n \rightarrow \infty}\frac{ \lvert \iota(\mathbf{v}_1) \rvert}{ \lvert \iota(\mathbf{u}) \rvert}  = i\sqrt{\lim_{n \rightarrow \infty}\frac{ \lvert \iota(\mathbf{v}_1) \rvert^2}{ \lvert \iota(\mathbf{u}) \rvert^2}}.
\]
Recall that $\lvert \iota(\mathbf{v}_1) \rvert ^2 = 6(\log n)^2 + O(\log n)$ and $\lvert \iota(\mathbf{u}) \rvert ^2 = 2(\log n)^2 + O(\log n)$. We now apply Lemma~\ref{lem:limit-helper} with $g(n) = \lvert \iota(\mathbf{v}_1) \rvert ^2$, $f(n) = 6(\log n)^2$, $h(n) = \lvert \iota(\mathbf{u}) \rvert ^2$, and $s(n) = 2(\log n)^2$. The hypotheses of Lemma~\ref{lem:limit-helper} are satisfied because 
\[
    g(x) = \lvert \iota(\mathbf{v}_1) \rvert ^2 = 6(\log n)^2 + O(\log n) = f(x) + O(\log n),
\]
and the last expression is $f(x) + o(f(x))$. Upon applying Lemma \ref{lem:limit-helper}, we obtain that:
\[
    \lim_{n \rightarrow \infty}\frac{ \lvert \iota(\mathbf{v}_1) \rvert^2}{ \lvert \iota(\mathbf{u}) \rvert^2} = 3
\]
so $\lim_{p \rightarrow \infty}p_L = i\sqrt{3}$.
\epf

\tpoint{Lemma} {
\label{lem:limit-helper}
Suppose we are given two functions $g(x),h(x)$. Suppose that $g(x) = f(x) + o(f(x))$ for some function $f$. Suppose that $h(x) = s(x) + o(s(x))$ some function $s$. Then
\[
    \lim_{x\rightarrow \infty} \frac{g(x)}{h(x)}
\]
exists if and only if
\[
    \lim_{x\rightarrow \infty} \frac{f(x)}{s(x)}
\]
exists. If both limits exist, then they are equal. 
}

\bpf
The assumption that $g(x) = f(x) + o(f(x))$ implies that $g$ and $f$ are asymptotically equivalent approaching infinity. Similarly, the assumption $h(x) = s(x) + o(s(x))$ implies that $h$ and $s$ are asymptotically equivalent approaching infinity.
\epf

\bpoint{The point $\frac{1}{7} + \frac{4i\sqrt{3}}{7}\in \overline{\Omega(D_4, 2,1)}$}

Again, we make use of Theorem~\ref{thm:stender} which describes the units of some pure quartic fields. Let $n$ be a positive integer such that $n^4 - 1$ is squarefree. Write $\omega = (n^4 - 1)^{1/4}$ and set $L_n = \QQ(\omega)$. 
By Theorem~\ref{thm:stender}, the units
\[
\left\{\iota\left(\frac{\omega+n}{\omega-n}\right), \iota\left(\frac{1}{\omega-n}\right)\right\}
\]
form a basis of $E_{L_n}$.

\tpoint{Proposition} { \em
\label{prop:third-point}
The set of limit points of $\Omega(D_4,2,1)$ contains $\frac{1}{7} + \frac{4i\sqrt{3}}{7}$.
}

\bpf
Choose $n$ as above. Now, let $\mathbf{u} = \frac{\omega+n}{\omega-n}$ and $\mathbf{v} = \frac{1}{\omega-n}$. A short calculation shows that
\begin{equation}
\label{eqn:u-point3}
\Log(\iota(\mathbf{u})) = (\log \lvert \omega + n\rvert-\log\lvert \omega - n\rvert, \log\lvert \omega - n \rvert -\log\lvert \omega + n\rvert,0)
\end{equation}
\begin{equation}
\label{eqn:v-point3}
\Log(\iota(\mathbf{v})) = -(\log\lvert\omega - n\rvert, \log\lvert\omega + n\rvert, 2\log\lvert i\omega + n\rvert).
\end{equation}
Set $\mathbf{u}_1 \coloneqq \mathbf{u}\mathbf{v}^{-1}$ and $\mathbf{v}_1 = \mathbf{v}^{-1}$. Clearly, $\{\iota(\mathbf{u}_1),\iota(\mathbf{v}_1)\}$ forms a basis of $E_{L_n}$. From Equation~\ref{eqn:u-point3} and Equation~\ref{eqn:v-point3}, we obtain:
\begin{equation}
\label{eqn:u1-point3}
\Log(\iota(\mathbf{u}_1))=  (\log\lvert\omega + n\rvert, \log\lvert\omega - n\rvert, 2\log\lvert i\omega + n\rvert)
\end{equation}
\begin{equation}
\label{eqn:v1-point3}
\Log(\iota(\mathbf{v}_1)) =  (\log\lvert\omega - n\rvert, \log\lvert\omega + n\rvert, 2\log\lvert i\omega + n\rvert).
\end{equation}
Now, the two vectors $\iota(\mathbf{u}_1)$ and $\iota(\mathbf{v}_1)$ have the same length. Let $0 < \theta \leq \pi$ denote the angle between the two vectors. If $\pi/3 \leq \theta \leq \pi/2$, then $\iota(\mathbf{u}_1)$ and $\iota(\mathbf{v}_1)$ form a reduced basis of the unit lattice of $L_n$. 

\noindent \textbf{Step 1: Calculating $\theta$.}
It is easy to see that
\begin{equation}
\label{eqn:bound1-point3}
\log\lvert\omega + n\rvert, \log \lvert i\omega + n \rvert = \log n + O(1).
\end{equation}
\begin{equation*}
\omega - n = (n^4 - 1)^{1/4} - n = n\Big((1 - \frac{1}{n^4})^{1/4} - 1\Big),
\end{equation*}
and the Taylor series expansion of $f(x) = (1 - x)^{1/4}$ around $0$ is
\[
    1 + \sum_{k = 1}^{\infty} \frac{f^k(0)}{k!}x^k = 1  - \frac{x}{4}+  \sum_{k = 2}^{\infty} \frac{(-1)^k(\frac{1}{4}\cdot \frac{3}{4} \cdot\dots \cdot \frac{-1 + 4k}{4})}{k!}x^k.
\]
Substituting, we get
\[
\omega - n = n\Big(1  - \frac{1}{4n^4}+  \sum_{k = 2}^{\infty} \frac{(-1)^k(\frac{1}{4}\cdot \frac{3}{4} \cdot\dots \cdot \frac{-1 + 4k}{4})}{n^{4k}k!} - 1 \Big) = \frac{1}{4n^3} + O(n^{-7}).
\]
which implies that
\begin{equation}
\label{eqn:bound2-point3}
\log\lvert\omega - n\rvert = -3\log(n) + O(1).
\end{equation}
From Equation~\ref{eqn:u1-point3} and Equation~\ref{eqn:v1-point3}, we obtain:
\begin{equation}
\label{eqn:cos-point3}
\cos\theta= \frac{\iota(\mathbf{u}_1) \cdot \iota(\mathbf{v}_1)}{\lvert \iota(\mathbf{u}_1) \rvert \lvert \iota(\mathbf{v}_1) \rvert}=\frac{2\log\lvert\omega + n\rvert\log\lvert\omega - n\rvert + 4\log^2\lvert i\omega + n\rvert}{\log^2\lvert\omega + n\rvert+ \log^2\lvert\omega - n\rvert + 4\log^2\lvert i\omega + n\rvert}.
\end{equation}

Now, substitute Equation~\ref{eqn:bound1-point3} and Equation~\ref{eqn:bound2-point3} into Equation~\ref{eqn:cos-point3}:
\[
\lim_{n\rightarrow \infty} \cos\theta  = \lim_{n\rightarrow \infty}\frac{-6\log(n)\log(n) + 4\log^2(n)}{\log^2(n) + 9\log^2(n) + 4\log^2(n)} = -1/7
\]
we obtain that if there are infinitely many such $n$, then $\lim_{n\to \infty} \cos \theta = \frac{-1}{7}$.

\noindent \textbf{Step 2: Showing that there are infinitely many fields.} To complete the proof, it suffices to show that there are infinitely many positive integers $n$ such that $n^4 - 1$ is squarefree. Because $n^4 - 1 = (n-1)(n+1)(n^2 + 1)$ has irreducible factors of degree $\leq 3$, Theorem 1.1 of \cite{browning-booker} implies that there are infinitely many such $n$.

Thus, we obtain a limit point on the unit circle (lower boundary) and have $\lim_{n \rightarrow \infty }p_L = \frac{1}{7} + \frac{4i\sqrt{3}}{7}$.
\epf
\bpoint{The point $\frac{1}{2}(1 + i\sqrt{3})\in \overline{\Omega(D_4, 2,1)}$}

To prove Proposition~\ref{prop:second-point}, we use a construction given in Section 2 of \cite{lps}, which we begin by expositing; we encourage the reader to refer to \cite{lps} for more details. 

For $n\in\ZZ$, we define
\[
    f_{n}(x) \coloneqq  x^4-nx^3+(1-n)x^2-2x+1.   
\]

When $f_{n}$ is irreducible, let $L_{n} = \QQ[x]/f_{n}(x)$ be the associated quartic field, and let $\rho \in L_{n}$ denote a root of $f$. If $\lvert n \rvert \geq 3$, then $n^2 - 4$ is a positive integer that is not a square; let $K_{n} = \QQ(\sqrt{n^2 - 4})$ be the associated real quadratic field. Set $\epsilon = \frac{1}{2}(-n + \sqrt{n^2 - 4})$; it is a unit of the ring of integers of $K_{n}$.

\tpoint{Theorem 2.6 of \cite{lps}} { \em
\label{thm:lps}
    If $n \leq -5$ and $(4n+17)(n^2-4)$ is squarefree, then the following four statements hold:
    \begin{enumerate}
        \item $f_{n}$ is irreducible and $L_{n}$ is a $D_4$-quartic field with signature $(2,1)$;
        \item $L_{n}$ contains the real quadratic subfield $K_{n}$;
        \item $\ZZ[\rho]$ is the maximal order of $L_{n}$;
        \item and $\{\iota(\epsilon), \iota(\rho)\}$ is a basis for the unit lattice of $L_{n}$.
    \end{enumerate}
}

We give a sketch of the proof of point $(4)$ in the theorem above. The polynomial $f_n$ is engineered so that $\rho$ and $\epsilon$ are units; the key point is to check that they span the unit lattice $E_{L_n}$. To do this, the authors use a bound of Nakamula \cite{nk} showing that the regulator of an order is bounded below by a certain function of the discriminant. This implies that $\rho$ and $\epsilon$ actually span the unit lattice. 

\tpoint{Proposition} { \em
\label{prop:second-point}
The set of limit points of $\Omega(D_4,2,1)$ contains $\frac{1}{2}(1 + i\sqrt{3})$.
}

\bpf
The key idea is to use Theorem~\ref{thm:lps} to construct a sequence of $D_4$-quartic fields $\{L_{n}\}$ such that $p_{L_{n}} \rightarrow \frac{1}{2}(1 + i\sqrt{3})$.

Let $n \leq -5$ be any integer. Suppose that $(4n+17)(n^2-4)$ is squarefree. Note that there exist infinitely many integers $n \leq -5$ for which $(4n+17)(n^2-4)$ is squarefree; this follows from a straightforward squarefree sieve, such as  Theorem 1 of \cite{erdos}. Now apply Theorem~\ref{thm:lps}.

\noindent \textbf{Step 1: Rewriting the basis.}
We have a factorization:
\begin{equation}
\label{eqn:factorization}
    f_{n} = (x^2 + \epsilon x + \epsilon)(x^2 - \overline{\epsilon} x + \overline{\epsilon})
\end{equation}
where $\epsilon = \frac{1}{2}(-n + \sqrt{n^2 - 4})$ and $\overline{\epsilon} = \frac{1}{2}(-n - \sqrt{n^2 - 4})$. Without loss of generality, we may suppose that $\rho$ is a root of $x^2 + \epsilon x + \epsilon$. Letting $\sigma$ be the nontrivial element of $\Gal(L_{n}/K_{n})$, Equation~\ref{eqn:factorization} shows that:
\[
    N_{L_{n}/K_{n}}(\rho) = \rho \sigma(\rho) = \epsilon.
\]
Therefore, because $\{\iota(\rho),\iota(\epsilon)\}$ is a basis for the unit lattice of $L_{n}$, so is $\{\iota(\rho), \iota(\sigma(\rho))\}$.

\noindent \noindent \textbf{Step 2: Computing $p_{L_{n}}$.} Now, the unit lattice $E_{L_{n}}$ is a rank $2$ lattice with an order two automorphism $\sigma$ which is not multiplication by $\pm 1$ (see Lemma~\ref{lem:nontriv-aut}). Moreover,  $\{\iota(\rho), \iota(\sigma(\rho)) = \sigma(\iota(\rho))\}$ is a basis for $E_L$; therefore, either:
\begin{enumerate}
    \item  $p_L$ lies on the lower boundary of $\cF$ and $\{\iota(\rho), \sigma(\iota(\rho))\}$ must be the two shortest vectors of $E_L$;
    \item or $p_L$ lies on the right boundary of $\cF$ and $\{\iota(\rho), \iota(N_{L/K}(\rho))\}$ is a short basis of $E_L$.
\end{enumerate}

Let $0 < \theta \leq \pi/2$ denote the angle between $\iota(\rho)$ and $\sigma(\iota(\rho))$. If $\theta > \frac{\pi}{3}$, then we must be in case $(1)$, i.e.  $p_L$ lies on the lower boundary of $\cF$. We now compute $\cos \theta$. 

Write $\Log(\iota(\rho)) = (a_1,a_2,-a_1-a_2)$ and $\Log(\sigma(\iota(\rho))) = (a_2,a_1,-a_1-a_2)$ for real numbers $a_1,a_2$. Computing the angle between $\Log(\iota(\rho))$ and $\Log(\sigma(\iota(\rho)))$, we obtain that: 
\begin{equation}
\label{eqn:cos}
 \cos \theta = \frac{ 2a_1a_2 + (a_1 + a_2)^2 }{ a_1^2 + a_2^2 + (a_1 + a_2)^2 }.
\end{equation}

Because $\rho$ is a root of  $x^2 + \epsilon x + \epsilon$, we apply the quadratic formula to obtain: 
\begin{equation}
\label{eqn:a1} 
 a_1  =  \log \bigg \lvert \frac{-\epsilon - \sqrt{\epsilon^2 - 4\epsilon}}{2} \bigg \rvert = \log \lvert \epsilon\rvert + O(1).  
\end{equation}

Similarly, we obtain an estimate for $a_2$:
\begin{equation}
\label{eqn:a2}
a_2 = \log \bigg \lvert \frac{-\epsilon + \sqrt{\epsilon^2  +4\epsilon}}{2} \bigg \rvert  = \log \lvert -\epsilon + \sqrt{\epsilon^2 + 4\epsilon} \rvert - \log(2)
\end{equation}

To estimate the first term in Equation~\ref{eqn:a2}, expand to get:
\begin{equation}
\label{eqn:a2-2}
-\epsilon + \sqrt{\epsilon^2 + 4\epsilon} = \epsilon\bigg(-1 + \sqrt{1 + \frac{4}{\epsilon}}\bigg). 
\end{equation}
Now the second term in Equation~\ref{eqn:a2-2} is bounded; we have $\lvert (-1 + \sqrt{1 + \frac{4}{\epsilon}})\rvert \asymp 2\epsilon^{-1}$, so combining this estimate with Equation~\ref{eqn:a2} and Equation~\ref{eqn:a2-2} implies that
\begin{equation}
\label{eqn:a2-bound}
    \lvert a_2 \rvert = O(1).
\end{equation}
Notice that as $n \rightarrow -\infty$, we have $\epsilon \rightarrow \infty$. Now, substituting Equation~\ref{eqn:a1} and Equation~\ref{eqn:a2-bound} into Equation~\ref{eqn:cos}, we obtain:
\begin{equation}
\label{eqn:cos-calc}
\lim_{n \rightarrow -\infty} \cos \theta = \frac{1}{2}
\end{equation}
as $n$ ranges along integers such that $(4n+17)(n^2-4)$ is squarefree. To proceed we'll need the following lemma. 

\tpoint{Lemma} {\em \label{angle-analysis}
    If $n$ is sufficiently large and $(4n+17)(n^2-4)$ is squarefree, then $\theta > \pi/3$.
}

\bpf 
We know that $\cos(\theta)$ approaches $\frac{1}{2}$ as $n \to -\infty$, but whether it approaches from the left or the right determines which boundary component the shapes lie on in the modular curve. To be precise, approaching from the left means that $p_{L_n}$ lies on the unit circle and $\theta \in (\pi/3, \pi/2)$. In this case the point lies on the lower boundary of the modular curve, we are in case $(1)$ and the basis is reduced. Approaching from the right means we are in case $(2)$, $\theta\in (0, \pi/3)$ and the basis is not reduced. In this case, the element $TST^{-1}\in \text{SL}_2(\mathbb{Z})$-where $S$ is inversion through the unit circle and $T$ is translation by one unit to the right-maps the point onto the right boundary of the modular curve. In fact, this element maps the entire arc, with $\theta\in (0, \pi/3)$, onto the vertical line with the real part $\frac{1}{2}$.
We show that we are in the former case.

By the computations above, we have:

\begin{equation}
    \cos(\theta)   = \frac{2a_1a_2 + (a_1 + a_2)^2}{a_1^2 + a_2^2 + (a_1 + a_2)^2} 
                   = \frac{a_1^2\left(1 + \frac{2a_2}{a_1} + \frac{a_2^2}{a_1^2}\right)}{2a_1^2\left(1 + \frac{a_2}{a_1} + \frac{a_2^2}{a_1^2}\right)} 
                   = \frac{1}{2} \left(\frac{1 + \frac{4a_2}{a_1} + \frac{a_2^2}{a_1^2}}{1 + \frac{a_2}{a_1} + \frac{a_2^2}{a_1^2}}\right).
\end{equation}
Now, since $a_1 \to \infty$ and $a_2 = O(1)$, we have $\frac{a_2}{a_1}\rightarrow 0$. Letting $x=\frac{a_2}{a_1}$ and expanding the resulting rational function yields:
\[
\frac{1 + 4x + x^2}{1 + x + x^2} = 1 + 3x - 3x^2 + O(x^4),
\]
and hence:
\[
\cos(\theta) = \frac{1}{2} \left(1 + 3x - 3x^2 +  O(x^4)\right).
\]

Since $x = \frac{a_2}{a_1} < 0$ for all sufficiently large $n$, we conclude that $\cos(\theta)$ approaches $\frac{1}{2}$ from the left. 
\epf 

By the lemma above, $\iota(\rho)$ and $\sigma(\iota(\rho))$ are a reduced basis of $E_L$. Therefore:
\begin{equation}
\label{eqn:theta}
p_{L_{n}} = \cos \theta + i\sin \theta.
\end{equation}

Now, substituting Equation~\ref{eqn:cos-calc} into Equation~\ref{eqn:theta}, we obtain:
\[
    \lim_{n \rightarrow -\infty} p_{L_n} \rightarrow \frac{1}{2} + \frac{i\sqrt{3}}{2}
\]
as $n$ ranges along integers such that $(4n+17)(n^2-4)$ is squarefree.
\epf


\bpoint{The point $ i\sqrt{3}\in \overline{\Omega(D_4, 2,1)}$}

For sake of completeness, we include a second proof that the set of limit points of $\Omega(D_4, 2,1)$ contains $i\sqrt{3}$. This comes from analyzing a different construction of explicitly parametrized units due to Nakamula \cite{nk}, which we begin by outlining here; we encourage the reader to consult Section $3$ of \cite{nk} for further details. For $n\in\ZZ_{>0}$, we define
\[
    f_{n}(x) \coloneqq  x^4-nx^3+3x^2-nx+1.   
\]

\tpoint{Lemma (Lemma~2, \cite{nk})} { \em
\label{lem:nk}
    Suppose $n > 3$, $n^2 - 4$ is not a square, $n^2 - 4 \equiv 1 \pmod{4}$, $25-4n^2$ is not a square, and $25 - 4n^2 \equiv 1 \pmod{4}$ and $\QQ(\sqrt{n^2 - 4}) \neq \QQ(\sqrt{25-4n^2})$. Then $f_{n}$ is irreducible and $L_{n} \coloneqq \QQ[x]/f_n(x)$ is a $D_4$-quartic field with signature $(2,1)$ and $L_{n}$ contains the real quadratic subfield $\QQ(\sqrt{n^2-4})$.
}

\tpoint{Lemma (Lemma~3, \cite{nk})} {\em
\label{lem:infinitude-of-n}
There are infinitely many $n$ satisfying the assumptions of Lemma~\ref{lem:nk}.
}

Now, in situation of Lemma~\ref{lem:nk}, we may explicitly describe the units of $L_n$. Let $\rho \in L_n$ denote the root of $f_n$.

\tpoint{Theorem (Proposition~4, \cite{nk})} { \em
\label{thm:nk-2}
    Under the assumptions of Lemma~\ref{lem:nk}, a basis of $E_{L_n}$ is given by $\{\iota(\rho), \iota(\rho + \rho^{-1})\}$.
}

We give a sketch of the proof of the above theorem. Again, the polynomial $f_n$ is explicitly constructed so that $\rho$ and $\rho + \rho^{-1}$ are units. As in the proof of Theorem~\ref{thm:lps}, Nakamula proves a lower bound on the regulator of a number field in terms of the discriminant. This lower bound shows that the two units $\{\iota(\rho), \iota(\rho + \rho^{-1})\}$ actually form a basis of $E_L$.

\tpoint{Proposition} { \em
\label{prop:fourth-point}
The set of limit points of $\Omega(D_4,2,1)$ contains $i\sqrt{3}$.
}

\bpf
We will construct a sequence of $D_4$-quartic fields $\{L_{n}\}$ such that $p_{L_{n}} \rightarrow i\sqrt{3}$. Let $n$ denote an integer satisfying the assumptions of Lemma~\ref{lem:nk}. Note that Lemma~\ref{lem:infinitude-of-n} shows that there are infinitely many such $n$.

Because $f_n$ is a reciprocal polynomial, if $f_n(\alpha) = 0$ for some nonzero $\alpha \in \overline{\QQ}$, then $f_n(\alpha^{-1}) = 0$. Let $\epsilon$ denote a real root of $f_n$ with $\lvert \epsilon \rvert \geq 1$, and let $\epsilon^{-1}$ be its inverse. Let $\eta, \eta^{-1}$ denote the complex roots of $f_n$. As usual, let $\rho \in L_n$ denote the root of $f_n$.

We begin by computing $\Log(\iota(\rho))$ and $\Log(\iota(\rho+\rho^{-1}))$. Order the real embeddings $\sigma_1,
\sigma_2$ of $L_n$ so that $\sigma_1(\rho) = \epsilon$. Then:
\begin{equation*}
     \Log(\iota(\rho)) = (\log \lvert \epsilon \rvert, -\log \lvert \epsilon \rvert, 2\log \lvert \eta \rvert).
\end{equation*}
Because the sum of the coordinates of $\Log(\iota(\rho))$ must be $0$, we must have $\log \lvert \eta \rvert = 0$, so we may write
\begin{equation}
\label{eqn:embedding-1}
     \Log(\iota(\rho)) = (\log \lvert \epsilon \rvert, -\log \lvert \epsilon \rvert, 0).
\end{equation}
Similarly, we have:
\begin{equation}
\label{eqn:embedding-2}
     \Log(\iota(\rho+\rho^{-1})) = \big(\log\lvert \epsilon+\epsilon^{-1} \rvert, \log\lvert \epsilon+\epsilon^{-1} \rvert, -2\log\lvert \epsilon+\epsilon^{-1} \rvert\big).
\end{equation}

A computation shows that $ \Log(\iota(\rho)) \cdot \Log(\iota(\rho+\rho^{-1})) = 0$, which implies that $E_L$ is an orthogonal lattice, so $p_{L_{n}}$ lies on the left boundary of $\cF$. Therefore, if $ \frac{\rvert\iota(\rho + \rho^{-1})\lvert}{\lvert\iota(\rho)\rvert} \geq 1$, then
\begin{equation}
\label{eqn:pln}
 p_{L_n} = i\frac{\lvert\iota(\rho + \rho^{-1})\rvert}{\lvert \iota(\rho)\rvert}.
\end{equation}

Combining Equation~\ref{eqn:embedding-1} and Equation~\ref{eqn:embedding-2}, we obtain:
\begin{equation}
\label{eqn:ratio-pt4}
\frac{\lvert\iota(\rho + \rho^{-1})\rvert}{\lvert\iota(\rho)\rvert} = \Bigg \lvert \frac{\sqrt{6}\log \lvert \epsilon + \epsilon^{-1} \rvert}{\sqrt{2}\log \lvert \epsilon \rvert} \Bigg \rvert =  \Bigg \lvert\frac{\sqrt{3}\log \lvert \epsilon + \epsilon^{-1}\rvert}{\log \lvert \epsilon \rvert} \Bigg \rvert.  
\end{equation}

An explicit computation of the discriminant of $f_n$ shows that $\lim_{n \rightarrow \infty} \lvert \Disc(f_n) \rvert = \infty$. Now, 
\[
    \Disc(f_n) = (\epsilon - \epsilon^{-1})^2(\epsilon - \eta)^2(\epsilon - \eta^{-1})^2(\epsilon^{-1} - \eta)^2(\epsilon^{-1} - \eta^{-1})^2(\eta - \eta^{-1})^2.
\]
In the computation in Equation~\ref{eqn:embedding-1}, we showed that $\lvert \eta \rvert = 1$ because $\log \lvert \eta \rvert = 0$. So, because $\lim_{n \rightarrow \infty} \lvert \Disc(f_n) \rvert = \infty$, we must have $\lim_{n\rightarrow \infty}\lvert \epsilon \rvert = \infty$. Hence, $\log\lvert\epsilon + \epsilon^{-1}\rvert = \log \lvert \epsilon \rvert + o(1)$. Therefore:
\begin{equation}
\label{eqn:final-eqn-pt4}
  \lim_{n\to\infty} \frac{\sqrt{3}\log\lvert \epsilon + \epsilon^{-1}\rvert} {\log \lvert \epsilon \rvert}= \sqrt{3}.  
\end{equation}
Substituting Equation~\ref{eqn:final-eqn-pt4} and Equation~\ref{eqn:ratio-pt4} into Equation~\ref{eqn:pln}, we obtain that
\[
    \lim_{n \rightarrow \infty} p_{L_n} = i\sqrt{3}.
\]
\epf

\end{document}